\documentclass[11pt]{article}
\usepackage{hyperref}
\usepackage{amsmath,amssymb}
\usepackage{ntheorem}
\usepackage{mathtools}
\usepackage{epsf}
\usepackage{theorem}
\usepackage{indentfirst}
\usepackage{latexsym}
\usepackage[dvips]{graphicx}
\usepackage{flafter}
\usepackage{caption}
\usepackage{textcomp}
\usepackage{hhline}
\usepackage{bigstrut,bigdelim}
\usepackage{rotating}
\usepackage{multicol,multirow}
\usepackage{makeidx}
\usepackage{epsfig}
\usepackage{fancyhdr}
\usepackage{booktabs}
\usepackage{mathrsfs}
\usepackage[all]{xy}
\usepackage{xcolor}
\usepackage{microtype}
\newif\ifger
\gerfalse
\newtheorem{definition}{Definition}[section]
\newtheorem{theorem}{Theorem}[section]
\newtheorem{lemma}[definition]{Lemma}
\newtheorem{corollary}[theorem]{Corollary}
\newtheorem{remark}{Remark}[section]

\newtheorem{proposition}[definition]{Proposition}
\newtheorem{example}{Example} [section]
\newtheorem*{funding}{Funding}
\newtheorem*{availability of data and material}{Availability of Data and Material}
\newtheorem*{conflict of interest}{Conflict of interest}
\newcommand{\Tr}{{\rm Tr}}
\newcommand{\gf}{ {{\mathbb F}} }
\begin{document}
\baselineskip=19pt

\title{Algebraic Structure of Permutational Polynomials\\ over $\mathbb{F}_{q^n}$ \uppercase\expandafter{\romannumeral2}}

\author{ Pingzhi Yuan$^a$, Xuan Pang$^a${\footnote{Corresponding author. E-mail: pangxuan202503@163.com(X. Pang) }, Danyao Wu$^b$} \\
\small \it  $^a$School of Mathematical Sciences, South China Normal University, \\
\small \it  Guangzhou, 510631, P. R. China\\
\small \it  $^b$School of Computer Science and Technology, Dongguan University of Technology, \\
\small \it Dongguan, 523808, P. R. China\\}

\date{}
\maketitle

\begin{abstract}
It is well known that there exists a significant equivalence between the vector space $\mathbb{F}_{q}^n$ and the finite fields $\mathbb{F}_{q^n}$, and many scholars often view them as the same in most contexts. However, the precise connections between them still remain mysterious. In this paper, we first show their connections from an algebraic perspective, and then propose a more general algebraic framework theorem. Furthermore, as an application of this generalized algebraic structure, we give some classes of permutation polynomials over $\mathbb{F}_{q^2}$.

\medskip
\noindent{\bf MSC(2020):} 11T06; 11T55.

\medskip
\noindent{\bf Keywords:} Finite field; permutation polynomial; algebraic framework theorem; algebraic structure.
\end{abstract}

\section{Introduction}

Let $\mathbb{F}_q$ be the finite field with $q$ elements, where $q=p^k$ is a prime power, and let $\mathbb{F}_q[x]$ be the ring of polynomials in a single indeterminate $x$ over $\mathbb{F}_q$. A polynomial $f(x)\in \mathbb{F}_{q}[x]$ is called a {\it permutation polynomial} (PP) of $\mathbb{F}_{q}$ if the induced map $\alpha \mapsto f(\alpha)$ from $\mathbb{F}_{q}$ to itself is bijective. Such polynomials have unique algebraic properties such as invertibility, that is, there exists a unique compositional inverse polynomial $f^{-1}(x)$ that satisfies $f\circ f^{-1}(x)=f^{-1}\circ f(x)=x$ in the sense of reduction modulo $x^q-x$. PPs over finite fields are of great significance in cryptography \cite{cryp1-RSA,cryp2-Schwenk,cryp3-Singh}, coding theory \cite{code1-Ding13,code2-DZ14,code3-Harish,code4-Laigle-Chapuy} and combinatorial design theory \cite{cd-DY06}, etc., and are powerful and indispensable tools in these domains.  


 Let $n$ be a positive integer. For any $x\in\mathbb{F}_{q^n}$,  the trace function $\Tr_{q}^{q^n}(x)$ from $\mathbb{F}_{q^n}$ to $\mathbb{F}_q$  is defined by
$$\Tr_{q}^{q^n}(x)=x+x^q+\dots+x^{q^{n-1}}.$$
We simply denote $\Tr_{q}^{q^n}(x)$ as $\Tr(x)$ for convenience. Let $\mathbb{F}_q^n$ be the $n$-dimensional vector space over $\mathbb{F}_q$, and $\mathbb{F}_q[x_1,\cdots,x_n]$ the polynomial ring in $n$ indeterminates over $\mathbb{F}_q$. Similarly, a function $g\in\mathbb{F}_q[x_1,\cdots,x_n]$ can also be defined as a permutation of $\mathbb{F}_q^n$ in the same manner. For more details on permutations, the reader is referred to the book of Lidl and Niederreiter \cite{Lidl}.

Recalled that $\mathbb{F}_{q^n}$ can also be regarded as an $n$-dimensional vector space over $\mathbb{F}_q$, which implies that PPs over $\mathbb{F}_{q^n}$ have certain algebraic structures. In light of this, we have explored and showed some algebraic structures of PPs from the algebraic standpoint. And Theorem 1.2 in \cite{yuan2024algebra1} may be considered as our starting point.
In most recently, several papers have emerged that focus on constructing PPs via the vector space
$\mathbb{F}_{q}^n$ \cite{qu2025,xiongms2024}. Although we often identify the vector space $\mathbb{F}_{q}^n$ with the finite fields $\mathbb{F}_{q^n}$, the precise connection between them is actually not clear to us.
 Inspired by this, we hope to build bridges between them, uncovering their deeper connections. And in this way, we deepen our exploration within the framework of $\mathbb{F}_q$-vector space homomorphism, thereby deriving the following generalized result.

\begin{theorem*}
\label{AST}
Given a PP $f(x)\in \mathbb{F}_{q^n}[x]$. Let $\rho:\mathbb{F}_{q^n}\rightarrow \mathbb{F}_q^n$, $\eta:\mathbb{F}_{q}^n\rightarrow\mathbb{F}_{q^n}$ be $\mathbb{F}_q$-vector space homomorphisms, and let $g$ be a map from $\mathbb{F}_q^n$ to $\mathbb{F}_q^n$. Then for any $x\in\mathbb{F}_{q^n}$ and  $(x_1,\dots,x_n)\in\mathbb{F}_q^n$, we have $$\rho(x)=(\Tr(v_1x),\dots,\Tr(v_nx)), \quad v_i\in\mathbb{F}_{q^n},~1\leq i\leq n,$$
$$\eta(x_1,\dots,x_n)=a_1x_1+\dots+a_nx_n, \quad a_i\in\mathbb{F}_{q^n},~ 1\leq i\leq n.$$
Moreover, $F(x)=\eta\circ g\circ \rho \circ f(x)$ is a PP of $\mathbb{F}_{q^n}$ if and only if the following two conditions hold:
\begin{enumerate}
\item[\rm(i)]Both $\{a_1,\dots, a_n\}$ and $\{v_1,\dots, v_n\}$ are bases of $\mathbb{F}_{q^n}$ over $\mathbb{F}_q$ ;
\item[\rm(ii)] $g\in \mathbb{F}_q[x_1,\dots, x_n]$ permutes $\mathbb{F}_{q}^n$.
\end{enumerate}
\end{theorem*}

Among the extensive studies of PPs, those of the form $x^rh(x^{q-1})$ over $\mathbb{F}_{q^2}$, with $r$ being a positive integer, have attracted remarkable interest over the past few decades, especially for PPs with fewer terms, such as permutation binomials  and permutation trinomials \cite{binomials-1,trinomial-1,trinomial-2,trinomial-3,trinomials,binomials-2,binomials-3,trinomials-4}.
 Notably, PPs with more terms have been considered in recent years.
Ding and Zieve \cite{zieve2023LondonMS} identified all permutation quadrinomials over $\mathbb{F}_{q^2}$ with the shape $x^rh(x^{q-1})$, under the condition that $r\equiv p^i+1\pmod {q+1}$ and $h(x)$ with special degrees. And in the most recent paper, the two authors \cite{ding_zieve2025} further investigated PPs over $\mathbb{F}_{q^2}$ of this form, and constructed two large classes of PPs, most of which have five terms. Zhang et al. \cite{zhenglijing2024FFA} and Kousar et al. \cite{xiongms2024} dedicated to the construction of permutation pentanomilas, and gave some PPs over $\mathbb{F}_{q^2}$ as well. For more results, see \cite{Shen2025,XuFFA2018,zhenglijiFFA2023}.  Therefore, we intend to construct several classes of PPs over $\mathbb{F}_{q^2}$, by utilizing the newly established algebraic framework theorem of PPs in our paper.

\begin{theorem}
Let $q$ be a prime power and $m,n\in \mathbb{N}$ be positive integers. Pick distinct $\alpha, \beta\in\mathbb{F}_{q^2}^*$ satisfying $\alpha^{q+1}=\beta^{q+1}=1$. Additionally, when $3\nmid q$, choose an element $\omega\in\mathbb{F}_{q^2}$ of order 3, and define $\varepsilon$ as follows: if $q\not\equiv2\pmod3$, let $\varepsilon\in\mathbb{F}_{q}^*$; if $q\equiv2\pmod3$, let $\varepsilon\in\mathbb{F}_{q}^*\cup\{\pm\omega, \pm\omega^2\}$. Under these conditions, the results presented in Table \ref{table1} hold.

\begin{table}[htbp]
  \centering
  \caption{Permutation polynomials over $\mathbb{F}_{q^2}$ in this paper}\label{table1}
  \begin{tabular}{lccl}
    \toprule
 {\rm No.} & $q$ & {\rm PPs} & {\rm N.S. Condition} \\
 \midrule
$1$ & & $(x+\alpha x^q)^m+\varepsilon(x+\beta x^q)^n$ & $\gcd(mn,q-1)=1$ {\rm and $\varepsilon^{q-1}\alpha^m\neq\beta^n$} \\
$2$ & $q\equiv 1\pmod3$ &   $(x+ x^q)^m+\varepsilon(x\pm\omega x^q)^n$  & $\gcd(mn,q-1)=1$ \\
$3$ &     &   $(x+\omega x^q)^m+\varepsilon(\omega x+ x^q)^n$  &   $\gcd(mn,q-1)=\gcd(3,m-2n)=1$ \\
$4$ &    &   $(x+\omega x^q)^m+\varepsilon(\omega x- x^q)^n$  &   $\gcd(mn,q-1)=1$ \\
$5$ & $q\equiv 2\pmod3$ &   $(x+ x^q)^m+\varepsilon( x+\omega x^q)^n$  & {\rm See Proposition \ref{PPq-2-1}}\\
$6$ &     &   $(x+ x^q)^m+\varepsilon( x-\omega x^q)^n$  & $\gcd(mn,q-1)=1$\\
$7$ &      &   $(x+\omega x^q)^m+\varepsilon(\omega x+ x^q)^n$  & {\rm See Proposition \ref{PPq-2-2}} \\
$8$ &       &   $(x+\omega x^q)^m+\varepsilon(\omega x-x^q)^n$  &  $\gcd(mn,q-1)=1$ \\
\bottomrule
  \end{tabular}
\end{table}
\end{theorem}

\begin{remark}
The N.S. condition in Table \ref{table1}, an abbreviation for {\it necessary and sufficient condition}, characterizes when the polynomial in column 3 is a PP over $\mathbb{F}_{q^2}$.
\end{remark}

The remainder of this paper is organized as follows. In Section 2,  we focus on uncovering the connections between the vector space $\mathbb{F}_{q}^n$ and the finite field $\mathbb{F}_{q^n}$. Building on these connections, we establish a generalized algebraic framework theorem and derive several useful corollaries. Then as an application of the algebraic framework theorem, Section 3 is devoted to the construction of PPs over $\mathbb{F}_{q^2}$ of the form $x^rh(x^{q-1})$. In Section 6, we conclude the paper.

\section{Algebraic structures of PPs}

Before proceeding to the proof of our algebraic framework theorem, (that is, Theorem \ref{AST}), we state the following known lemma and then give two technical lemmas that will be useful below.

\begin{lemma}\cite[Lemma 4.2]{yuan2024algebra1} \label{lem2.1}
Let $f(x)\in\mathbb{F}_{q^n}[x]$ be a map from $\mathbb{F}_{q^n}$ to $\mathbb{F}_q$. Then $f(x)$ is an $\mathbb{F}_q$-space homomorphism if and only if $f(x)=\Tr(vx)$ for some $v\in\mathbb{F}_{q^n}$.
\end{lemma}

\begin{lemma}\label{lem2.2}
Let $\rho:\mathbb{F}_{q^n}\rightarrow \mathbb{F}_q^n$ be an $\mathbb{F}_q$-vector space homomorphism, then there are $n$ elements $v_1,\dots,v_n \in \mathbb{F}_{q^n}$ such that
\begin{equation}\label{eq1}
\rho(x)=(\Tr(v_1x),\dots, \Tr(v_nx)),\quad \forall x \in \mathbb{F}_{q^n}.
\end{equation}
Moreover, $\rho$ is an $\mathbb{F}_q$-vector space isomorphism if and only if $\{v_1,\dots, v_n\}$ is a basis of $\mathbb{F}_{q^n}$ over $\mathbb{F}_q$, and when this occurs, its compositional inverse $\rho^{-1}:\mathbb{F}_{q}^n\to \mathbb{F}_{q^n}$ is given by:
  \begin{equation}\label{eq2}
  \rho^{-1}(x_1,\dots,x_n)=u_1x_1+\dots +u_nx_n,\quad \forall (x_1,\dots,x_n)\in\mathbb{F}_{q}^n,
  \end{equation}
where $\{u_1,\dots,u_n\}$ is the ordered dual basis of $\{v_1,\dots,v_n\}$.

\end{lemma}
\textbf{Proof.}\quad Suppose that $\rho(x)=(f_1(x),\cdots, f_n(x))$, where each $f_i(x)\in \mathbb{F}_q$. We define a set of mappings $\Delta=\{\sigma_1,\cdots, \sigma_n\}$ where
  $\sigma_i:\mathbb{F}_{q}^n\rightarrow\mathbb{F}_q$  is a map defined by $(x_1,\dots,x_n)\mapsto x_i$. Observe that for each $i$, the composition $\sigma_i\circ \rho$ yields an $\mathbb{F}_q$-vector space homomorphism from $\mathbb{F}_{q^n}$ to $\mathbb{F}_q$. Therefore, by using Lemma \ref{lem2.1},  we conclude that there exists an element $v_i\in\mathbb{F}_{q^n}$ such that $f_i(x)=\sigma_i\circ\rho(x)=\Tr(v_ix)$. This completes the proof of the first part.

Now we prove the sufficiency of the second part. Let $\{v_1,\dots, v_n\}$ be a basis of $\mathbb{F}_{q^n}$ over $\mathbb{F}_q$, and let $\{u_1,\dots, u_n\}$ be its dual basis with respect to the trace form. 
 Consequently, for each element $u_i$, we have $\rho(u_i)=e_i$, where
$e_i=(0,\dots,0,1,0,\dots,0)\in\mathbb{F}_{q}^n$ denotes the $i$-th standard unit vector in the $n$-dimensional vector space $\mathbb{F}_q^n$. Therefore, $\rho$ is an $\mathbb{F}_q$-vector space isomorphism since $\rho$ maps a basis $\{u_1,\dots, u_n\}$ of $\mathbb{F}_{q^n}$ to the standard basis $\{e_1,\dots, e_n\}$ of $\mathbb{F}_q^n$.

Conversely, let $\rho$ be an $\mathbb{F}_q$-vector space isomorphism. Assume that $\{v_1,\dots, v_n\}$ is not a basis of $\mathbb{F}_{q^n}$ over $\mathbb{F}_{q}$, and then there exists a set of maximal linearly independent subset within $\{v_1,\dots, v_n\}$. Without loss of generality, let $\{v_1,\dots,v_r\}$ with $r<n$ be such a maximal subset, which implies that any $v_i$ for $i=r+1,\dots,n$, can be written as a linear combination of the vectors in $\{v_1,\dots,v_r\}$ with coefficients from $\mathbb{F}_q$. Moreover, since $\mathbb{F}_{q^n}$ is $n$-dimensional over $\mathbb{F}_q$, it is always possible to extend  $\{v_1,\dots,v_r\}$ to be a basis of $\mathbb{F}_{q^n}$. Specially, there exist $\{v_{r+1}^*,\dots, v_n^*\}\in\mathbb{F}_{q^n}$ such that
$\{v_1,\cdots,v_r,v_{r+1}^*,\dots, v_n^*\}$
forms a basis of $\mathbb{F}_{q^n}$ over $\mathbb{F}_q$. Subsequently, we can take its ordered dual basis $\{u_1,\dots,u_n\}$. Now consider  the kernel of $\rho$. we have $$ker\rho=\{x:\Tr(v_ix)=0,1\leq i\leq n\}=<u_{r+1},\dots,u_n>\neq \{0\},$$
which contradicts the fact that $\rho$ is an $\mathbb{F}_q$-vector space isomorphism.

It remains to verify that $\rho^{-1}$ as defined in (\ref{eq2}) is indeed the compositional inverse of $\rho$. On the one hand, for any $x\in \mathbb{F}_{q^n}$, we have $x=\sum_{i=1}^n{u_ik_i}$ since $\{u_1,\dots,u_n\}$ is a basis of $\mathbb{F}_{q^n}$ over $\mathbb{F}_q$, where each $k_i\in\mathbb{F}_q$. Recall that $\{u_1,\dots,u_n\}$ and $\{v_1,\dots,v_n\}$ are a dual pair of ordered bases of $\mathbb{F}_{q^n}$ over $\mathbb{F}_q$, it follows that $Tr(v_jx)=k_i$ for $1\leq j\leq n$, and thus
$$\rho^{-1}\circ \rho(x)=\rho^{-1}(\Tr(v_1x),\dots,\Tr(v_nx))=u_1\Tr(v_1x)+\dots+u_n\Tr(v_nx)=x,$$
implying that $\rho^{-1}\circ \rho=id_{\mathbb{F}_{q^n}}$. On the other hand, it is easy to check that $\rho\circ \rho^{-1}=id_{\mathbb{F}_{q}^n}$.

This completes the proof of the lemma. $\hfill\square$

\begin{lemma}\label{lem2.3}
Let $\eta:\mathbb{F}_{q}^n\rightarrow \mathbb{F}_{q^n}$  be an $\mathbb{F}_q$-vector space homomorphism, then there exist $n$ elements  $a_1,\dots, a_n \in \mathbb{F}_{q^n}$ such that
\begin{equation}\label{eq3}
\eta(x_1,\dots,x_n)=a_1x_1+\dots+a_nx_n, \quad \forall (x_1,\dots, x_n)\in\mathbb{F}_{q}^n.
\end{equation}
Moreover, $\eta$ is an $\mathbb{F}_{q}$-vector space isomorphism if and only if $\{a_1,\dots,a_n\}$ is a basis of $\mathbb{F}_{q^n}$ over $\mathbb{F}_q$, and when this occurs, its compositional inverse $\eta^{-1}:\mathbb{F}_{q^n}\to \mathbb{F}_{q}^n$ is given by:
  \begin{equation}\label{eq4}
  \eta^{-1}(x)=(\Tr(b_1x),\dots, \Tr(b_nx)),\quad \forall x\in\mathbb{F}_{q^n},
  \end{equation}
where $\{b_1,\dots,b_n\}$ is the ordered dual basis of $\{a_1,\dots,a_n\}$.
\end{lemma}
\textbf{Proof.}\quad Let $\eta$ be an $\mathbb{F}_q$-vector space homomorphism. For any standard unit vector $e_i\in\mathbb{F}_q^n$, there exists a unique element $a_i\in \mathbb{F}_{q^n}$ such that $\eta(e_i)=a_i$. By the linearity property of $\eta$, for any $(x_1,\dots,x_n)\in \mathbb{F}_{q}^n$, we have
$$\eta(x_1,\dots,x_n)=\eta(x_1e_1+\dots+x_ne_n)=a_1x_1+\dots+a_nx_n.$$
Furthermore, $\eta$ is an $\mathbb{F}_q$-vector space isomorphism if and only if  the images of the standard basis vector $\{e_1,\dots, e_n\}$ form a basis of $\mathbb{F}_{q^n}$. In other word, $\{a_1,\dots, a_n\}$ is a basis of $\mathbb{F}_{q^n}$ over $\mathbb{F}_{q}$.

Finally, with a similar argument as in the proof of Lemma \ref{lem2.1}, we conclude that $\eta^{-1}$ defined in (\ref{eq4}) is indeed the compositional inverse of $\eta$.
$\hfill\square$

\bigskip

With these preparations complete, we are now ready to finish the proof of Theorem \ref{AST}.

\textbf{Proof of Theorem  \ref{AST}:}

\vspace{0.5em}

Proof.\quad Indeed, the first past of the theorem has already been proved. To complete the proof, it suffices to address the second part. Based on the given conditions, $f(x)$ is a PP of $\mathbb{F}_{q^n}$. We define $\rho: \mathbb{F}_{q^n}\to \mathbb{F}_{q}^n$, $\eta:\mathbb{F}_{q}^n\to\mathbb{F}_{q^n}$  to be $\mathbb{F}_q$-vector space homomorphisms, and the map $g:\mathbb{F}_{q}^n\to\mathbb{F}_{q}^n$ given by $$(x_1,\dots,x_n)\mapsto(g_1(x_1,\cdots,x_n),\dots,g_n(x_1,\dots,x_n)).$$
Clearly, $\mathbb{F}_{q^n}$, as an $n$-degree extension field of $\mathbb{F}_q$, contains exactly $q^n$ elements, and $\mathbb{F}_q^n$, which is the Cartesian product of $n$ copies of $\mathbb{F}_q$, also consists of $q^n$ elements. So, the following commutative diagram holds.

\vspace{0.5em}

\centerline{\xymatrix@C=3.3em@R=2.8em{
\mathbb{F}_{q^n} \ar[r]^{\textstyle f} \ar@/_2em/[rrrr]_{ F(x)} & \mathbb{F}_{q^n} \ar[r]^{ \textstyle \rho}  & \mathbb{F}_{q}^n \ar[r]^{\textstyle g}& \mathbb{F}_q^n \ar[r]^{\textstyle \eta} & \mathbb{F}_{q^n}
}}
\vspace{0.5em}
\noindent It follows that $F(x)=\eta\circ g\circ\rho\circ f(x)$ is a PP of $\mathbb{F}_{q^n}$ if and only of $\rho$, $\eta$ are $\mathbb{F}_q$-vector space isomorphisms, and $g$ permutes $\mathbb{F}_q^n$. And thus by applying Lemma \ref{lem2.2} and Lemma \ref{lem2.3}, we immediately complete the proof. $\hfill\square$

\bigskip
We now state the following two corollaries as immediate consequences of Theorem \ref{AST}.
Essentially these corollaries imply an important bidirectional correspondence, that is, for every permutation of $\mathbb{F}_{q^n}$, it is always possible to construct a  permutation of $\mathbb{F}_{q}^n$, and vice versa.

\begin{corollary}
Given $\rho:\mathbb{F}_{q^n}\to \mathbb{F}_q^n$ and $\eta:\mathbb{F}_q^n\to \mathbb{F}_{q^n}$ be two $\mathbb{F}_q$-vector space homomorphisms, and a mapping $g:\mathbb{F}_q^n\to\mathbb{F}_q^n$. Then the composition  $\eta\circ g\circ\rho$ induces a map from $\mathbb{F}_{q^n}$ to $\mathbb{F}_{q^n}$. Furthermore, if $g$ is a permutation over $\mathbb{F}_q^n$, and $\rho$, $\eta$ are $\mathbb{F}_q$-vector space isomorphisms, the mapping $\eta\circ g\circ\rho$ is a PP over $\mathbb{F}_{q^n}$.
\end{corollary}

\begin{corollary}
Given $\rho:\mathbb{F}_{q^n}\to \mathbb{F}_q^n$, $\eta:\mathbb{F}_q^n\to \mathbb{F}_{q^n}$ be two $\mathbb{F}_q$-vector space homomorphisms, and a mapping $f:\mathbb{F}_{q^n}\to\mathbb{F}_{q^n}$. Then the composition $\rho\circ f\circ\eta$ induces a map from $\mathbb{F}_{q}^n$ to $\mathbb{F}_{q}^n$. Furthermore, if $f$ is a PP of $\mathbb{F}_{q^n}$, and both $\rho$ and $\eta$ are $\mathbb{F}_q$-vector space isomorphisms, the mapping $\rho\circ f\circ\eta$ is a permutation of  $\mathbb{F}_{q}^n$.
\end{corollary}
\textbf{Proof.}\quad  To see this, we merely need to make the following slight modification based on the commutative diagram above.

\vspace{0.3em}
\centerline{\xymatrix@C=3.3em@R=2.8em{
\mathbb{F}_{q}^n \ar[r]^{\textstyle \eta}  & \mathbb{F}_{q^n} \ar[r]^{ \textstyle f}  & \mathbb{F}_{q^n} \ar[r]^{\textstyle \rho}& \mathbb{F}_q^n
}} $\hfill\square$

To illustrate the applicability of these two corollaries, we present a simple example below.
\begin{example}
Let $\{a_1,\dots, a_n\}$ and $\{v_1,\dots, v_n\}$ be a dual pair of ordered bases of $\mathbb{F}_{q^n}$ over $\mathbb{F}_q$, and let $\rho$ and $\eta$ be $\mathbb{F}_q$-vector space isomorphisms defined as (\ref{eq1}) and (\ref{eq3}), respectively. Then for the canonical PP  $f(x)=x\in\mathbb{F}_{q^n}[x]$, the composition
$$\rho\circ f\circ\eta(x)=(\Tr(v_1\sum_{i=1}^{n}{a_ix_i}),\dots,\Tr(v_n\sum_{i=1}^{n}{a_ix_i}) )=(x_1,\dots, x_n),$$
acts as the identity map on the vector space $\mathbb{F}_{q}^n$.
\end{example}

Further, we take into account the two fundamental polynomial rings, namely $\mathbb{F}_{q^n}[x]$ and $\mathbb{F}_q[x_1,\dots,x_n]$, defined over the finite field $\mathbb{F}_q$, and a profound algebraic property is then exhibited as follows.
For simplicity, let  $Map(A, B)$ denote the set of all mappings from $A$ to $B$.

\begin{proposition}
Let $\left (Map(\mathbb{F}_{q}^n, \mathbb{F}_{q}^n),+,\cdot\right)$ and $ \left(Map(\mathbb{F}_{q^n}, \mathbb{F}_{q^n}),+,\cdot\right)$ be $\mathbb{F}_q$-vector spaces, equipped with addition and scalar multiplication operations. Then $\psi$ given by
\begin{align*}
\psi: Map(\mathbb{F}_{q}^n, \mathbb{F}_{q}^n)&\to Map(\mathbb{F}_{q^n}, \mathbb{F}_{q^n})\\
g & \mapsto \rho^{-1}g\rho
\end{align*}
is an $\mathbb{F}_q$-vector space isomorphism, where $\rho:\mathbb{F}_{q^n}\to \mathbb{F}_{q}^n$ and $\rho^{-1}:\mathbb{F}_{q}^n\to \mathbb{F}_{q^n}$ are $\mathbb{F}_q$-vector isomorphisms defined as (\ref{eq1}) and (\ref{eq2}), respectively. Additionally, $\psi$ is also a groupoid isomorphism with respect to the composition.
\end{proposition}
\textbf{Proof.}\quad Obviously, the map $\psi$ is one-to-one. And it is easy to verify that for arbitrary $a_1,a_2\in\mathbb{F}_q$ and $g_1,g_2\in Map(\mathbb{F}_{q}^n, \mathbb{F}_{q}^n)$, the following property holds:
$$\psi(a_1g_1+a_2g_2)=a_1\psi(g_1)+a_2\psi(g_2).$$
Hence, $\psi$ is an $\mathbb{F}_q$-vector space isomorphism. Furthermore, it follows from $\rho^{-1}(g_1\circ g_2)\rho=(\rho^{-1}g_1\rho )\circ (\rho^{-1} g_2\rho)$ that
$$\psi(g_1\circ g_2)=\psi(g_1)\circ\psi(g_2),$$
which shows that $\psi$ respects the groupoid composition law. This finishes the proof.
$\hfill\square$

\bigskip
In the study of PPs, one with a simple structure is desired. Hence we define the mapping $g:\mathbb{F}_q^n\to \mathbb{F}_q^n$ with the following special form: $(x_1,\dots,x_n)\mapsto (h_1(x_1),\cdots,h_n(x_n))$, i.e., $g_i(x_1,\cdots,x_n)=h_i(x_i)\in\mathbb{F}_q[x_i]$,  for $1\leq i\leq n$. Then by Theorem \ref{AST}, we obtain the following statement, which improves our result in \cite{yuan2024algebra1}.

\begin{theorem}\label{AST2}
Given a PP $f(x)\in \mathbb{F}_{q^n}[x]$, and polynomials $h_i(x)\in\mathbb{F}_{q}[x]$, $1\leq i\leq n$. Then the polynomial
$$F(x)=a_1h_1(\Tr(v_1f(x)))+\dots+a_nh_n(\Tr(v_nf(x))),~~a_i,v_i\in\mathbb{F}_{q^n}, ~1\leq i\leq n$$
is a PP of $\mathbb{F}_{q^n}$ if and only if the following two conditions hold:
\begin{enumerate}

\item[\rm(i)]Both $\{a_1,\dots, a_n\}$ and $\{v_1,\dots v_n\}$ are bases of $\mathbb{F}_{q^n}$ over $\mathbb{F}_q$;
\item[\rm(ii)] $h_i(x)$, $1\leq i\leq n$, are PPs over $\mathbb{F}_q$.
\end{enumerate}
\end{theorem}
\textbf{Proof.}\quad For clarity, we show a clearer diagram below.

\vspace{0.5em}

\centerline{
  \xymatrix@C=6.6em@R=2.8em{
    \mathbb{F}_{q^n} \ar[r]^{f(x)}
    & \mathbb{F}_{q^n} \ar[r]^{\rho}_{\raisebox{-1em}{$\scriptstyle x \,\mapsto\, (\Tr(v_1x_1), \dots,\, \Tr(v_nx_n))$}}
    & \mathbb{F}_{q}^n \ar[r]^{g}_{\substack{ (x_1,\dots,x_n) \,\mapsto\, \\ (h_1(x_1),\ldots,h_n(x_n)) }}
    & \mathbb{F}_q^n \ar[r]^{\eta}_{\substack{ (x_1,\dots,x_n) \,\mapsto\, \\ a_1x_1 + \dots + a_nx_n }}
    & \mathbb{F}_{q^n}
  }
}   $\hfill\square$

\begin{remark}
 Theorem \ref{AST2} remains valid for the mapping $g:\mathbb{F}_{q}^n\to \mathbb{F}_{q}^n$ defined as
 $$g(x_1,\dots,x_n)= (h_1(x_1),h_2(x_2)+g_2(x_1),h_2(x_2)+g_3(x_1,x_2),\dots, h_n(x_n)+g_n(x_1,\dots,x_{n-1})),$$
 where each component satisfies $h_i(x_i)+g_i(x_1,\dots,x_{i-1})\in \mathbb{F}_q[x_1,\dots,x_i]$. In fact, the validity can be derived from the fact that the $i$-component depends soely on $x_i$.
\end{remark}

It is well known that the function $x^m$ is a PP of $\mathbb{F}_{q}$ if and only if $\gcd(m,q-1)=1$. Applying this elementary result, we then deduce a corollary as follows.
\begin{corollary}
If $f(x)\in\mathbb{F}_{q^n}[x]$ is a PP over $\mathbb{F}_{q^n}$ and $m_i\in \mathbb{N}$, $1\leq i\leq n$ are positive integers, then the polynomial
$$F(x)=a_1(\Tr(v_1f(x)))^{m_1}+\cdots+a_n(\Tr(v_1f(x)))^{m_n},~~a_i,v_i\in\mathbb{F}_{q^n}, ~1\leq i\leq n$$
is a PP over $\mathbb{F}_{q^n}$ if and only if $\gcd(m_1\dots m_n,q-1)=1$, and both $\{a_1,\dots, a_n\}$ and $\{v_1,\dots, v_n\}$ are bases of $\mathbb{F}_{q^n}$ over $\mathbb{F}_q$.

\end{corollary}

\section{PPs of the form $x^rh(x^{q-1})$ over $\mathbb{F}_{q^2}$}

 In this section, we would like to construct some PPs over $\mathbb{F}_{q^2}$.  We first introduce some basic tools that will underpin the analyses and proofs in the remainder of this paper.

%
%
%
%
%
%
%
%
%
%

\begin{lemma}\label{lem3.1}
Let $q$ be a prime power. Pick an order $3$ element $\omega\in\gf_{q^2}$ if $3\nmid q$. Then

(1) If $a\in\gf_{q^2}$ is an element with $a^q=a\alpha$, where $\alpha\in\mathbb{F}_{q^2}^*$, then $\Tr(ax)=a(x+\alpha x^q)$.

(2) If $3\nmid q$ and $a\in\gf_{q^2}$ is an element with $a^q=a\omega$, then $\Tr(ax)=a(x+\omega x^q)$.

(3) If $3\nmid q$ and $a\in\gf_{q^2}$ is an element with $a^q=-a\omega$,  then $\Tr(ax)=a(x-\omega x^q)$.

(4) If $q\equiv1\pmod{3}$ and $a\in\gf_{q^2}$ is an element with $a^q=a\omega^2$,  then $\Tr(a\omega x)=a(\omega x+ x^q)$.

(5) If $q\equiv2\pmod{3}$ and $a\in\gf_{q^2}$ is an element with $a^q=a\omega$,  then $\Tr(a\omega x)=a(\omega x+ x^q)$.

(6) If $q\equiv1\pmod{3}$ and $a\in\gf_{q^2}$ is an element with $a^q=-a\omega^2$,  then $\Tr(a\omega x)=a(\omega x- x^q)$.

(7) If $q\equiv2\pmod{3}$ and $a\in\gf_{q^2}$ is an element with $a^q=-a\omega$,  then $\Tr(a\omega x)=a(\omega x- x^q)$.

\end{lemma}

The Dimension Theorem guarantees that in a two-dimensional vector space, two linearly independent vectors span the space and thus form a basis. Hence by Theorem \ref{AST2}, we have the following corollary.
\begin{corollary}\label{cor3.1}
Let $q$ be a prime power, and let $a_1,a_2,b_1,b_2\in\mathbb{F}_{q^2}$, and $g_1(x),g_2(x)\in \mathbb{F}_q[x]$.  Then $f(x)=b_1g_1(\Tr(a_1x))+b_2g_2(\Tr(a_2x))$
 permutes $\mathbb{F}_{q^2}$ if and only if
 \begin{enumerate}
\item[\rm(i)]Both  $\{a_1,a_2\}$ and $\{b_1,b_2\}$ are $\mathbb{F}_q$-linearly independent;
\item[\rm(ii)] $g_1(x),g_2(x)$ are PPs over $\mathbb{F}_q$.
\end{enumerate}
\end{corollary}

\begin{remark}\label{linearly-inde}
To demonstrate the $\mathbb{F}_q$-linear independence of two elements $a_1,a_2\in \mathbb{F}_{q^2}$, it suffices to verify that the ratio $(a_2/a_1) \notin\mathbb{F}_q$, equivalently,  $(a_2/a_1)^{q-1}\neq 1$.

\end{remark}

In the subsequent propositions, we construct several classes of PPs over $\mathbb{F}_{q^2}$, serving as concrete applications of our proposed algebraic framework theorem.  We note, however, that the first result presented below is not new. Indeed, it was previously proved in our work via  an investigation of the compositional inverse of the polynomial \cite{yuan2025inverses}. Moreover, Ding and Zieve \cite{ding_zieve2023,zieve2013} obtained such kind of PP in the special case where $m=n$ as well.

\begin{proposition}\label{PP1}
Let $q$ be a prime power and $m,n$ be positive integers. Pick distinct $\alpha, \beta\in\mathbb{F}_{q^2}^*$ such that $\alpha^{q+1}=\beta^{q+1}=1$. Then $$f(x)=(x+\alpha x^q)^m+\varepsilon(x+\beta x^q)^n,$$
where $\varepsilon\in\mathbb{F}_{q^2}^*$, is a PP over $\mathbb{F}_{q^2}$ if and only if $\gcd(mn,q-1)=1$ and $\varepsilon^{q-1} \alpha^m\neq\beta^n$.
\end{proposition}
\textbf{Proof.}\quad As shown in part (1) of Lemma \ref{lem3.1}, we choose two elements $a_1, a_2\in\mathbb{F}_{q^2}$ satisfying $a_1^{q-1}=\alpha$ and $a_2^{q-1}=\beta$. Then we can see that
$$f(x)=(x+\alpha x^q)^m+\varepsilon(x+\beta x^q)^n=a_1^{-m}\Tr^m(a_1x)+\varepsilon a_2^{-n}\Tr^n(a_2x).$$

\noindent Consequently, by Corollary \ref{cor3.1}, $f(x)$ is a PP of $\mathbb{F}_{q^2}$ if and only if $\gcd(mn, q-1)=1$, and both the sets $\{a_1,a_2\}$ and $\{a_1^{-m}, \varepsilon a_2^{-n}\}$  are linearly independent over $\mathbb{F}_q$.
From Remark \ref{linearly-inde},
$\{a_1^{-m}, \varepsilon a_2^{-n}\}$   is $\mathbb{F}_q$- linearly independent if and only if
$(a_1^{-m}/\varepsilon a_2^{-n})^{q-1}\neq1$, which is equivalent to $\varepsilon^{q-1} \alpha^m\neq\beta^n$. On the other hand, $\{a_1,a_2\}$  must be $\mathbb{F}_q$- linearly independent since $\alpha\neq\beta$.
 This complete the proof.
$\hfill\square$

We hereby give some specific examples by applying the above proposition.
\begin{example}\label{example}
Let $q=p^k$ be an odd prime power. Pick an element $\alpha$ of order $q+1$ in $\mathbb{F}_{q^2}^*$, and let $\beta=-\alpha$. If  $\varepsilon=1$ and $m=n=3$, then we have
$$x^3+3\alpha^2x^{1+2q}$$
is a PP over $\mathbb{F}_{q^2}$ if and only if $\gcd(3,q-1)=1$ .
If $\varepsilon=1$ and $m=n=5$. We have
$$x^5+10\alpha^2x^{3+2q}+5\alpha^4x^{1+4q}$$
is a PP over $\mathbb{F}_{q^2}$ if and only if $\gcd(5,q-1)=1$ .
If $\varepsilon=1$ and $m=n=7$. We have
$$x^7+21\alpha^2x^{5+2q}+35\alpha^4x^{3+4q}+7\alpha^6x^{1+6q}$$
is a PP over $\mathbb{F}_{q^2}$ if and only if $\gcd(7,q-1)=1$ .
If $\varepsilon=1$ and  $m=n=Q+R+S$, where $Q,R,S\in\{p^i:i\geq0\}$. We have
$$x^{Q+R+S}+\alpha^{R+S}x^{Q+q(R+S)}+\alpha^{Q+S}x^{R+q(Q+S)}+\alpha^{Q+R}x^{S+q(Q+R)}$$
is a PP over $\mathbb{F}_{q^2}$ if and only if $\gcd(Q+R+S,q-1)=1$.
\end{example}

\begin{remark}
%

Lappano \cite[Theorem 3.3.1]{Lappano} proved that binomial $f(x)=ax^3+x^{1+2q}\in\mathbb{F}_q[x]$ with $a\neq0$ is a PP of $\mathbb{F}_{q^2}$ if and only if (i) $a=1$ and $q\equiv1\pmod4$; (ii) $a=1/3$ and $q\equiv-1\pmod6$; or (iii) $a=-1/3$ and $q\equiv-1\pmod{12}$.
 Notably, part (ii) and (iii) are in accord with our findings, provided that we choose $\alpha\in\mu_{q+1}$  such that $\alpha^2=1$ for part (ii) and $\alpha^2=-1$ for part (iii). And the existence of $\alpha$ satisfying $\alpha^2=-1$ can be justified as follows. Since  $\alpha^{q+1}=1$, it follows that $(\alpha^{\frac{q+1}{2}})^2=1$. If $q\equiv-1\pmod{12}$, then $\frac{q+1}{2}$ is even, and it is possible to find an element $\alpha\in\mu_{q+1}$ such that $\alpha^2=-1$. In addition, it is worth noting that several classes of PPs involving trinomials or quadrinomials are shown in Example \ref{example} as well, and we believe that these examples likely contain some new families of PPs, which further enrich the known results in this field.

\end{remark}

\begin{proposition}
Let $q\equiv1\pmod{3}$ be a prime power and $m,n$ be positive integers. Pick an order $3$ element $\omega\in\gf_{q}$. Then
$$f(x)=(x+ x^q)^m+\varepsilon( x\pm\omega x^q)^n,$$
where $\varepsilon\in\mathbb{F}_q^*$, is a PP over $\mathbb{F}_{q^2}$ if and only if $\gcd(mn,q-1)=1$.
\end{proposition}
\textbf{Proof.}\quad  We begin by proving for the case where the signs of the two terms in $f(x)$ are identical. From Lemma \ref{lem3.1}, part (2), we select an element $a\in\mathbb{F}_{q^2}$ satisfying $a^{q-1}= \omega$ and then the polynomial
$$f(x)=(x+ x^q)^m+\varepsilon( x+\omega x^q)^n=\Tr^m(x)+\varepsilon a^{-n}\Tr^n(ax).$$

\noindent By Corollary \ref{cor3.1}, it suffices to prove that both $\{1,a\}$ and $\{1,\varepsilon a^{-n}\}$ form bases of $\mathbb{F}_{q^2}$ over $\mathbb{F}_q$, provided that $\gcd(mn,q-1)=1$. Assume now this condition holds. It follows from $q\equiv1\pmod3$ that $3\nmid n$, and
$$(\varepsilon a^{-n})^{q-1}=\omega^{-n}\neq1,$$
since $\varepsilon\in\mathbb{F}_q^*$ and $\omega$ is a primitive 3rd root of unity. By Remark \ref{linearly-inde}, we are done.

 For another case in which the terms have opposite signs, we choose another element $\hat{a}$ to replace $a$, where $\hat{a}\in\mathbb{F}_{q^2}$ satisfies   $\hat{a}^{q-1}=-\omega$. Then, for any positive integer $n$, we have
$$ (\varepsilon \hat{a}^{-n})^{q-1}=-\omega^{-n}\neq1.$$
Therefore, $\{1, \hat{a}\}$ and $\{1, \varepsilon \hat{a}^{-n}\}$ are necessarily two bases of $\mathbb{F}_{q^2}$ over $\mathbb{F}_q$.  Once more, by Corollary \ref{cor3.1}, we conclude the proof. $\hfill\square$

%
%

\begin{remark}
In the case where the signs of the two terms in $f(x)$ are opposite,  we take $\mathbb{F}_{q^2}$ as a finite field of odd characteristic by default. Otherwise, the concept of ``opposite signs" becomes trivial since $+1=-1$ in the fields of even characteristic.
\end{remark}

\begin{proposition}\label{PP3}
Let $q\equiv1\pmod{3}$ be a prime power and $m, n$ be positive integers. Pick an order $3$ element $\omega\in\gf_{q}$. Then
$$f(x)=(x+ \omega x^q)^m+\varepsilon(\omega x\pm x^q)^n,$$
where $\varepsilon\in\mathbb{F}_q^*$, is a PP over $\mathbb{F}_{q^2}$ if and only if $\gcd(mn,q-1)=\gcd(3,m-2n)=1$.
\end{proposition}
\textbf{Proof.} \quad  Let us first consider the case of $f(x)$ with identical signs.
Applying parts (2) and (4) of Lemma \ref{lem3.1}, we can express $f(x)$ as:
$$f(x)=(x+\omega x^q)^m+\varepsilon(\omega x+ x^q)^n=a_1^{-m}\Tr^m(a_1x)+\varepsilon a_2^{-n}\Tr^n(a_2\omega x),$$
where $a_1, a_2$ are two elements of $\mathbb{F}_{q^2}$ with $a_1^{q-1}=\omega$, $a_2^{q-1}=\omega^2$. Now we assert that $\{a_1,a_2\omega\}$ and $\{a_1^{-n},\varepsilon a_2^{-n}\}$ are $\mathbb{F}_{q}$-linearly independent, provided that $\gcd(mn,q-1)=\gcd(3,m-2n)=1$. This follows from
$$(a_2\omega/a_1)^{q-1}=\omega, \quad (\varepsilon a_2^{-n}/a_1^{-m})^{q-1}=\omega^{m-2n}\neq1.$$
Consequently, by Corollary \ref{cor3.1}, we arrive at the desired result.

 For the remaining case, consider another element $\hat{a}_2\in\mathbb{F}_{q^2}$ satisfying $\hat{a}_2^{q-1}=-\omega^2$, along with previously-defined $a_1$, we have
 $$f(x)=(x+\omega x^q)^m+\varepsilon(\omega x- x^q)^n=a_1^{-m}\Tr^m(a_1x)+\varepsilon \hat{a}_2^{-n}\Tr^n(\hat{a}_2\omega x).$$
Given $\gcd(mn,q-1)=1$, we obtain the following equations:
 $$(\hat{a}_2/a_1)^{q-1}=-\omega, \quad (\varepsilon \hat{a}_2^{-n}/a_1^{-m})^{q-1}=-\omega^{m-2n}=(-\omega)^{m-2n}.$$
We now claim that $(-\omega)^{m-2n}\neq1$. Here, since $q\equiv1\pmod3$ is an odd prime power,  the element $-\omega\in\mathbb{F}_{q^2}$ is a primitive sixth root of unity. From the condition $\gcd(mn,q-1)=1$, we can deduce that $m-2n$ is odd.  Because any integer divisible by 6 must be even, it is clear that $m-2n$ cannot be divisible by 6. As a result, we have $(-\omega)^{m-2n}\neq1$. This complete the proof. $\hfill\square$

\begin{remark}\label{PP-independent}
In fact, the proof of Proposition \ref{PP3} shows that for the negative sign case, $f(x)$ is a PP over $\mathbb{F}_{q^2}$ if and only if $\gcd(mn,q-1)=1$, regardless of the condition $\gcd(3,m-2n)=1$.
\end{remark}

\begin{proposition}\label{PPq-2-1}
Let $q\equiv2\pmod{3}$ be a prime power and $m,n$ be positive integers. Pick an order $3$ element $\omega\in\gf_{q^2}$, and let
$$f(x)=(x+  x^q)^m+\varepsilon (x+ \omega x^q)^n,$$
where $\varepsilon\in\mathbb{F}_q^*\cup\{\pm\omega,\pm\omega^2\}.$
Then we have the following:
\begin{enumerate}
\item[\rm(i)] If $\varepsilon\in\mathbb{F}_q^* $, then $f(x)$ permutes $\mathbb{F}_{q^2}$ if and only if $(mn,q-1)=1$ and $3 \nmid n$;
\item[\rm(ii)] If $\varepsilon=\pm \omega$, then $f(x)$ permutes $\mathbb{F}_{q^2}$ if and only if $(mn,q-1)=1$ and $n\not\equiv1\pmod3$;
\item[\rm(iii)] If $\varepsilon=\pm \omega^2$, then $f(x)$ permutes $\mathbb{F}_{q^2}$ if and only if $(mn,q-1)=1$ and $n\not\equiv2\pmod3$.
\end{enumerate}
\end{proposition}
\textbf{Proof.} \quad
Case (i) is a direct consequence of Proposition \ref{PP1}, but for the sake of completeness, we provide its analysis below as well. Taking $a\in\mathbb{F}_{q^2}$ with  $a^{q-1}=\omega$,  then the polynomial can be expressed as
 $$f(x)=(x+ x^q)^m+\varepsilon ( x+\omega x^q)^n=\Tr^m(x)+\varepsilon a^{-n}\Tr^n(a x).$$
It is known that $f(x)$ is a PP of $\mathbb{F}_{q^2}$ if and only if $\gcd(mn,q-1)=1$, and both $\{1,a\}$ and $\{1,\varepsilon a^{-n}\}$ are bases of $\mathbb{F}_{q^2}$ over $\mathbb{F}_q$. Undoubtedly, the set $\{1,a\}$ forms a basis of $\mathbb{F}_{q^2}$ over $\mathbb{F}_q$.
Our goal is to determine necessary and sufficient conditions for the $\mathbb{F}_q$-linearly independence of the set $\{1,\varepsilon  a^{-n}\}$, i.e, consider the non-equality $(\varepsilon a^{-n})^{q-1}\neq1$. We achieve this through a classified discussion over $\varepsilon$.
If $\varepsilon\in\mathbb{F}_{q}^*$, then
$$(\varepsilon a^{-n})^{q-1}=\omega ^{-n}\neq1$$ if and only if $3\nmid n$.
 If $\varepsilon=\pm\omega$, then
$$(\varepsilon a^{-n})^{q-1}=(\pm\omega a^{-n})^{q-1}=\omega ^{1-n}\neq1$$ if and only if $n\not\equiv 1\pmod3$.
 If $\varepsilon=\pm\omega^2$, then
$$(\varepsilon a^{-n})^{q-1}=(\pm\omega^2 a^{-n})^{q-1}=\omega ^{2-n}\neq1$$ if and only if $n\not\equiv 2\pmod3$.
Above all, we arrive at the conclusion. $\hfill\square$

\begin{proposition}\label{PPq-2-1-negative}
Let $q\equiv2\pmod{3}$ be a prime power and $m,n$ be positive integers. Pick an order $3$ element $\omega\in\gf_{q^2}$, and let
$$f(x)=(x+  x^q)^m+\varepsilon (x- \omega x^q)^n,$$
where $\varepsilon\in\mathbb{F}_q^*\cup\{\pm\omega,\pm\omega^2\}.$ Then $f(x)$ is a PP over $\mathbb{F}_{q^2}$ if and only if $\gcd(mn,q-1)=1$.
\end{proposition}
\textbf{Proof.} \quad From part (3) of Lemma \ref{lem3.1}, we choose an element $a\in\mathbb{F}_{q^2}$ such that $a^{q-1}=-\omega$. Then $f(x)$ can be represented in the form $f(x)=\Tr^m(x)+\varepsilon a^{-n}\Tr(ax)$. Notice that for any $\varepsilon\in\mathbb{F}_{q}^*\cup\{\pm\omega,\pm\omega^2\}$, we have
$$(\varepsilon a^{-n})^{q-1}=\varepsilon^{q-1}(-\omega)^{-n}=-\varepsilon^{q-1}\omega^{-n}\neq1$$
The inequality follows from $\omega$ is an element of order 3, and consequently none of the elements in the set  $\{-\omega^{-n},-\omega^{1-n},-\omega^{2-n}\}$ equals 1. Thus by Corollary \ref{cor3.1}, this proposition follows immediately.

\begin{proposition}\label{PPq-2-2}
Let $q\equiv2\pmod{3}$ be a prime power and $m,n$ be positive integers. Pick an order $3$ element $\omega\in\gf_{q^2}$, and let
$$f(x)=(x+ \omega x^q)^m+\varepsilon (\omega x+ x^q)^n,$$
where $\varepsilon\in\mathbb{F}_{q}^*\cup\{\pm\omega,\pm\omega^2\}.$
Then we have the following:
\begin{enumerate}
\item[\rm(i)] If $\varepsilon\in\mathbb{F}_q^* $, then $f(x)$ permutes $\mathbb{F}_{q^2}$ if and only if $(mn,q-1)=1$ and $3\nmid (m-n)$;
\item[\rm(ii)] If $\varepsilon=\pm \omega$, then $f(x)$ permutes $\mathbb{F}_{q^2}$ if and only if $(mn,q-1)=1$ and  $3\nmid (1+m-n)$;
\item[\rm(iii)] If $\varepsilon=\pm \omega^2$, then $f(x)$ permutes $\mathbb{F}_{q^2}$ if and only if $(mn,q-1)=1$ and $3\nmid (2+m-n)$.
\end{enumerate}
\end{proposition}
\textbf{Proof.} \quad
Let $a_1,a_2\in\mathbb{F}_{q^2}$ be elements satisfying $a_1^{q-1}=\omega$ and $a_2^{q-1}=\omega$.  Then $f(x)$ is rewritten in the form of $f(x)=a_1^{-m}\Tr^m(a_1x)+\varepsilon a_2^{-n}\Tr^n(a_2\omega x)$. Next by an argument analogous to Proposition \ref{PPq-2-1}, this statement can be derived.
 $\hfill\square$

\begin{proposition}\label{PPq-2-2-negative}
Let $q\equiv2\pmod{3}$ be a prime power and $m,n$ be positive integers. Pick an order $3$ element $\omega\in\gf_{q^2}$, and let
$$f(x)=(x+ \omega x^q)^m+\varepsilon (\omega x- x^q)^n,$$
where $\varepsilon\in\mathbb{F}_{q}^*\cup\{\pm\omega,\pm\omega^2\}.$
$f(x)$ is a PP over $\mathbb{F}_{q^n}$ if and only if $\gcd(mn,q-1)=1$.
\end{proposition}
\textbf{Proof.} \quad
Choose two element $a_1,a_2\in\mathbb{F}_{q^2}$ satisfying $a_1^{q-1}=\omega$ and $a_2^{q-1}=-\omega$. Following a discussion similar to that in Proposition \ref{PPq-2-1-negative}, we obtain the validity of the stated proposition.
$\hfill\square$

\begin{remark}
 In \cite{ding_zieve2025},  Ding and  Zieve constructed some classes of PPs over $\mathbb{F}_{q^2}$ with the shape $x^rh(x^{q-1})$. Hereby we aim to use our approach to provide an alternative explanation of the permutation pentanomials appearing in their Theorem 1.7.

{\bf Case 1:} $q\equiv1\pmod{3}$.

For Theorem 1.1 in case $z = 1$,
induces the same function on $\gf_{q^2}$ as does
$$(X+\omega X^q)^{Q+R+S}-\omega(\omega X+X^q)^{Q+R+S},$$
which by inspection equals $\beta f(X)$. Now by computation, we have
$$(X+\omega X^q)^{Q+R+S}-\omega(\omega X+X^q)^{Q+R+S}=$$
$$a^{-Q-R-S}\left(\Tr(aX)\right)^{Q+R+S}-\omega b^{-Q-R-S}\left(\Tr(b\omega X)\right)^{Q+R+S},$$
where $a$ and $b$ are two elements of $\gf_{q^2}$ with
$a^q=a\omega$, $b^q=b\omega^2.$

For Theorem 1.3 in case $z = 1$, we have
$$\beta f(X)=(X+\omega X^q)^{Q+qR+S}-\omega(\omega X+X^q)^{Q+qR+S}.$$
 Now by computation, we have
$$(X+\omega X^q)^{Q+qR+S}-\omega(\omega X+X^q)^{Q+qR+S}=$$
$$a^{-Q-qR-S}\left(\Tr(aX)\right)^{Q+qR+S}-\omega b^{-Q-qR-S}\left(\Tr(b\omega X)\right)^{Q+qR+S},$$
where $a$ and $b$ are two elements of $\gf_{q^2}$ with
$a^q=a\omega$,  $b^q=b\omega^2.$

For Theorem 1.1 in case $z = 2$, we have
$$\beta f(X)=(\omega X+ X^q)^{Q+R+S}-\omega( X+\omega X^q)^{Q+R+S}.$$
 Now  through calculation, we have
$$(\omega X+ X^q)^{Q+R+S}-\omega( X+\omega X^q)^{Q+R+S}=$$
$$a^{-Q-R-S}\left(\Tr(a\omega X)\right)^{Q+R+S}-b^{-Q-R-S}\left(\Tr(b X)\right)^{Q+R+S},$$
where $a$ and $b$ are two elements of $\gf_{q^2}$ with
$a^q=a\omega^2$, $ b^q=b\omega.$

For Theorem 1.3 in case $z = 2$, we have
$$\beta f(X)=(\omega X+ X^q)^{Q+qR+S}-\omega( X+\omega X^q)^{Q+qR+S}.$$
 Now by computation, we have
$$(\omega X+ X^q)^{Q+qR+S}-\omega( X+\omega X^q)^{Q+qR+S}=$$
$$a^{-Q-qR-S}\left(\Tr(a\omega X)\right)^{Q+qR+S}-\omega b^{-Q-qR-S}\left(\Tr(b X)\right)^{Q+qR+S},$$
where $a$ and $b$ are two elements of $\gf_{q^2}$ with
$a^q=a\omega^2$, $ b^q=b\omega.$

{\bf Case 2:} $q\equiv2\pmod{3}$.

For Theorem 1.1 in case $z = 1$,
induces the same function on $\gf_{q^2}$ as does
$$(X+\omega X^q)^{Q+R+S}-\omega (\omega X+X^q)^{Q+R+S},$$
 by simple computation, we have
 $$ (X+\omega X^q)^{Q+R+S}-\omega (\omega X+X^q)^{Q+R+S}=$$
$$a^{-Q-R-S}(\Tr(aX))^{Q+R+S}-\omega b^{-Q-R-S}(\Tr(b\omega X))^{Q+R+S},$$
where $a$ and $b$ are two elements of $\gf_{q^2}$ with
$a^q=a\omega$,  $b^q=b\omega.$

For Theorem 1.1 in case $z = 2$, we have
$$(\omega X+ X^q)^{Q+R+S}-\omega (X+\omega X^q)^{Q+R+S}.$$
Now by computation, we have
 $$ (\omega X+ X^q)^{Q+R+S}-\omega (X+\omega X^q)^{Q+R+S}=$$
$$a^{-Q-R-S}(\Tr(a\omega X))^{Q+R+S}-\omega b^{-Q-R-S}(\Tr(bX))^{Q+R+S},$$
where $a$ and $b$ are two elements of $\gf_{q^2}$ with
$a^q=a\omega$,  $b^q=b\omega.$

As shown in Table \ref{table1}, these pentanomials given in \cite{ding_zieve2025}  indeed are PPs. Furthermore, for $q\equiv2\pmod3$, we easily show that adding $\alpha\in\mu_{q+1}$ to both two terms of the above polynomials does not affect on their permutation property.
 For example,
 $$(\omega X+\alpha X^q)^{Q+R+S}-\omega(X-\alpha \omega X^q)^{Q+R+S}$$
permutes $\mathbb{F}_{q^2}$ if and only if $\gcd(Q+R+S,q-1)=1$.
\end{remark}
\section{Conclusion and further work}
In this paper, we clearly explained the equivalence relationship between the vector space $\mathbb{F}_{q}^n$ and the finite fields $\mathbb{F}_{q^n}$ from an algebraic standpoint. Based on their connections, we provided a more generalized algebraic framework theorem than Theorem 1.2 in \cite{yuan2024algebra1}. In addition, we constructed several families of permutation polynomials of the form $x^rh(x^{q+1})$ over $\mathbb{F}_{q^2}$ listed in Table  \ref{table1}, by leveraging the generalized algebraic framework theorem. And we believe that more permutation polynomials could be obtained by our approach.

For further work, exploring planar functions and almost perfect nonlinear (APN) functions by using the isomorphism between $Map(\mathbb{F}_{q}^n, \mathbb{F}_{q}^n)$ and $Map(\mathbb{F}_{q^n}, \mathbb{F}_{q^n})$ may be a fascinating work. Additionally, while our current work only gave an alternative interpretation of Ding and Zieve's results in \cite{ding_zieve2025}, the fact is that more known permutation polynomials might be explained in a similar way.


\bigskip
\section*{Acknowledgements}

We are sincerely grateful to  Michael E. Zieve for his timely attention and for pointing out some errors in first version of our manuscript. These errors, if left uncorrected, could lead to inaccuracies in our work. Thanks to his comments, we were able to correct them in time and improve the overall quality and credibility of the manuscript.
\textit{ \begin{funding}{\rm Pingzhi Yuan was supported by the National Natural Science Foundation of China (Grant No. 12171163) and Guangdong Basic and Applied Basic Research Foundation  (Grant No. 2024A1515010589).
}
\end{funding}}
%


\section*{Declarations}
\begin{conflict of interest} {\rm There is no conflict of interest.}
\end{conflict of interest}

\end{document}